\documentclass{amsart}
 \usepackage{amsmath}
 \usepackage{amssymb}
 \usepackage{amscd}
 \usepackage[mathscr]{eucal}
 \usepackage{epsfig}
\usepackage{graphics} 
\usepackage{url}

 \newtheorem{thm}{Theorem}[section]
 
 \newtheorem{remark}[thm]{Remark}
 \newtheorem{example}[thm]{Example}
 \newtheorem{cor}[thm]{Corollary}
 \newtheorem{prop}[thm]{Proposition}
 \newtheorem{defn}[thm]{Definition}

\newcommand{\R}{\mathbb{R}}
\newcommand{\C}{\mathbb{C}}

\newcommand{\miss}{MISS}

\begin{document}
 
 \title[Spherical bundles]{Enumerative properties of triangulations of spherical bundles over $S^1$}
\author{Jacob Chestnut}
\address{Evans Hall, University of California, Berkeley, CA, 94720}
\email{ jacob\_chestnut@berkeley.edu}
\author{Jenya Sapir}
\address{Dept. of Mathematics, U. of Chicago, Chicago, IL 60637}
\email{ jmsapir@uchicago.edu}
\author{Ed Swartz}
\address{Malott Hall, Cornell Univerisity, Ithaca, NY 15853}
\email{ebs22@cornell.edu}

\begin{abstract}
We give a complete characterization of all possible pairs $(f_0,f_1)$, where $f_0$ is the number of vertices and $f_1$ is the number of edges, of  any  triangulation of an $S^k$-bundle  over $S^1.$  The main point is that K\"uhnel's triangulations of $S^{2k+1} \times S^1$ and the nonorientable $S^{2k}$-bundle over $S^1$ are unique among all triangulations of $(n-1)$-dimensional homology manifolds with $2n+1$ vertices, first Betti number nonzero, and whose orientation double cover has vanishing second Betti number.  
\end{abstract}

\maketitle

The basic enumerative invariant of any triangulation of an $(n-1)$-dimensional compact manifold is its $f$-vector, $(f_0,\dots,f_{n-1}),$ where $f_i$ is the number of $i$-dimensional simplices. There are very few manifolds for which  a complete description of all  $f$-vectors is known.  Indeed, in dimension five and above there are {\it no} manifolds where this problem has been solved.  The $f$-vectors of compact surfaces were determined by Ringel \cite{Ri}, and Jungerman and Ringel \cite{JR}.  In addition to $S^4,$ the three-manifolds $S^3, \R P^3, S^2 \times S^1$ and the nonorientable $S^2$ -bundle over $S^1$ were covered by Walkup \cite{Wal}. All possible $f$-vectors of $\C P^2, (S^2 \times S^2) \# (S^2 \times S^2), S^3 \times S^1$ and K3-surfaces were determined in \cite{Sw10}.  Our Theorem \ref{s3 twist s1} gives all  $f$-vectors of the nonorientable $S^3$ bundle over $S^1.$  

While a thorough understanding of $f_2$ (and higher) remains elusive in dimensions above four, all pairs $(f_0, f_1)$ are known for several higher dimensional manifolds.  Brehm and K\"uhnel proved that the minimum number of vertices for a PL-triangulation of $S^{2k+1} \times S^1$ is $4k+7$, while the minimum number of vertices for a PL-triangulation of the nonorientable $S^{2k}$-bundle over $S^1$ is $4k+5$ \cite{BKu}.  This was extended in \cite{Sw10} to all  triangulations for the same collection of spherical bundles over $S^1.$   In addition, all possible combinations for the number of vertices and edges were given for $S^{2k+1} \times S^1$ \cite[Proposition 5.4]{Sw10}.  Theorem \ref{edges+} provides a complete characterization of all possible pairs $(f_0, f_1)$ for all spherical bundles over $S^1.$

Minimal triangulations of $S^{2k+1} \times S^1$ and the nonorientable $S^{2k}$-bundles over $S^1$ were originally found by K\"uhnel \cite{Ku2}.  In 
Section \ref{unique} we prove that any triangulation of an $S^k$ bundle over $S^1$ with $2k+5$ vertices is combinatorially isomorphic to one of K\"uhnel's minimal triangulations.  The proof of Theorem \ref{edges+}  consists of combining this uniqueness result with the constructive methods in Sections \ref{algorithm} and \ref{edges}.

As will become apparent,  all of our constructions yield combinatorial manifolds.  Hence, our results also apply verbatim when restricted to this smaller class of triangulations.   Immediately after we wrote this paper we discovered the arXive preprint of Bhaskar Bagchi and Basudeb Datta, ``The lower bound theorem and minimal triangulations of sphere bundles over the circle,"  arXiv:math.GT/0610829, which has results very similar (but not identical) to ours \cite{BD}.    Bagchi and Datta consider the category of manifold triangulations and, in addition to constructing minimal  triangulations of spherical bundles over $S^1,$ prove that any non-simply connected $(n-1)$-dimensional  manifold with $2n+1$ vertices is isomorphic to one of K\"uhnel's minimal triangulations.  

 \section{Preliminaries}

Throughout, $\Delta$ is a connected, pure, $(n-1)$-dimensional simplicial complex with $m$ vertices and vertex set $V  = \{v_1, \dots, v_m\}.$ A simplicial complex is {\it pure} if all of its facets (maximal faces)  have the same dimension.  In addition, we will always  assume that $n \ge 4.$  The {\it geometric realization} of $\Delta, |\Delta|,$ is the union in $\R^m$ over all faces $\{v_{i_1}, \dots, v_{i_j}\}$ of $\Delta$ of the convex hull of $\{e_{i_1},\dots,e_{i_j}\},$ where $\{e_1, \dots, e_m\}$ is the standard basis of $\R^m.$ We say $\Delta$ is {\it homeomorphic} to another space whenever $|\Delta|$ is.  A {\it triangulation} of a topological space $M$ is any simplicial complex $\Delta$ such that $\Delta$ is homeomorphic to $M.$  

The {\it $f$-vector} of $\Delta$ is $( f_0, \dots, f_{n-1}),$ where $f_i$ is the number of $i$-dimensional faces in $\Delta.$ Sometimes it is convenient to set $f_{-1} = 1$, corresponding to the empty set.   The {\it face polynomial} of $\Delta$ is 

$$f_\Delta(x) = f_{-1} x^n + f_0 x^{n-1} + \dots + f_{n-2} x + f_{n-1}.$$

The {\it $h$-vector} of $\Delta$ is $(h_0, \dots, h_n)$ and is defined so that the corresponding $h$-polynomial, $h_\Delta(x) = h_0 x^n + h_1 x^{n-1} + \dots +h_{n-1} x +  h_n,$ satisfies $h_\Delta(x+1) = f_\Delta(x).$  Equivalently,

\begin{equation} \label{f by h}
  h_i = \sum^i_{j=0} (-1)^{i-j} \binom{n-j}{n-i} f_{j-1}.
\end{equation} 
 
Each $f_i$ is a  {\it nonnegative} linear combination of  $h_0, \dots, h_{i+1}.$  Specifically,
\begin{equation} \label{h by f}
  f_{i-1} = \sum^i_{j=0} \binom{n-j}{n-i} h_j.
\end{equation}

Evidently the $f$-vector and $h$-vector encode the same information. One of the  advantages of the $h$-vector  is that the linear equalities satisfied by triangulations of manifolds have a very simple form.

\begin{thm}\cite{Ke}
    If $\Delta$ is a triangulation of a closed manifold,  then
        \begin{equation} \label{klee}
        h_{n-i} - h_i =  (-1)^i \binom{n}{i} (\chi (\Delta) - \chi (S^{n-1})).
       \end{equation}
    \end{thm}
\noindent In fact, Klee's formula holds for the more general class of {\it semi-Eulerian} complexes.

When $\Delta$ is homeomorphic to a manifold, Klee's formula (\ref{klee}) allows us to specify the $f$-vector of $\Delta$ using only $h_0,\dots, h_{\lfloor n/2 \rfloor}.$  This is one of the motivations behind introducing the $g$-vector.  For $i \le \lfloor n/2 \rfloor$  define
$$g_i = h_i - h_{i-1}.$$
\noindent  In view of (\ref{klee}), the $f$-vector of a triangulation of a manifold is determined by its $g$-vector $(g_0,\dots,g_{\lfloor n/2 \rfloor}).$  As we will see, 
$$g_2 = h_2 - h_1 = f_1 - n f_0 + \binom{n+1}{2}$$ plays a special role.

A {\it stacked polytope} is the following inductively defined class of polytopes.   The simplex is a stacked polytope and any polytope obtained from a stacked polytope by adding a pyramid to a facet is a stacked polytope.  Stacked polytopes are simplicial and the boundary of a stacked polytope is a {\it stacked sphere}.   From a purely combinatorial point of view, a stacked sphere is obtained by beginning with the boundary of a simplex and then repeatedly subdividing facets, i.e.  replacing a facet $\{v_{i_1}, \dots, v_{i_n}\}$ with $n$ facets $\{v^\prime, v_{i_2}, \dots, v_{i_n}\}, \{v_{i_2}, v^\prime, v_{i_3}, \dots, v_{i_n}\}, \dots, \{v_{i_1}, \dots, v_{i_{n-1}}, v^\prime\},$ where $v^\prime$ is a new vertex. 

One method for constructing triangulations of spherical bundles over $S^1$ is to start with a triangulation $\Delta$ of $S^n,$ identify two facets and then remove the identified facet. We say the resulting space is obtained from $\Delta$ by {\it handle addition}.   As long as there is no path of length less than three between  each pair of identified vertices in $\Delta$, the resulting space will be a simplicial complex homeomorphic to an $S^{n-1}$-bundle over $S^1.$   As there are, up to homeomorphism, only two such bundles, $S^{n-1} \times S^1$ and a nonorientable space \cite{Ste}, the topological type of the quotient space is determined by the orientation of the identification.  When the original sphere is a stacked sphere we call such a space an {\it identified stacked sphere} or {\it ISS}.  The importance of ISS's is demonstrated by the following.  Here, $\beta_i$ is the $i^{th}$ Betti number with respect to rational coefficients.

\begin{thm} \label{MISS->SS} \cite[Theorem 4.30]{Sw10}  Suppose $\Delta$ is a triangulation of an oriented $(n-1)$ homology manifold with $\beta_1 \neq 0, \beta_2 = 0$   and $n \ge 5.$  Then $g_2 \ge \beta_1 \binom{n+1}{2}.$  Furthermore, if $g_2 = \beta_1 \binom{n+1}{2},$ then $\Delta$ is an ISS.    

\end{thm}

\begin{cor} \label{cor}
If $\Delta$ is a triangulation of an $S^{n-2}$-bundle over $S^1$ with $n \ge 5,$ then $g_2 \ge  \binom{n+1}{2}.$  Furthermore, if $g_2 =  \binom{n+1}{2},$ then $\Delta$ is an ISS.
\end{cor}

\begin{proof}
If $\Delta$ is orientable, then the above theorem applies, so assume that $\Delta$ is not orientable.  Let $\Delta^\prime$ be the induced triangulation on $S^{n-2} \times S^1,$ the orientation double cover of $\Delta.$  Direct computation and  the fact that $f_i(\Delta^\prime) = 2 f_i(\Delta)$ for all $i \ge 0$ implies that \cite[Proposition 4.2]{Sw10}
$$g_2(\Delta) = \frac{g_2(\Delta^\prime) + \binom{n+1}{2}}{2}.$$
Since $g_2(\Delta^\prime) \ge \binom{n+1}{2},$ so is $g_2(\Delta).$  If $g_2(\Delta) = \binom{n+1}{2},$ then $g_2(\Delta^\prime) = \binom{n+1}{2}$ which implies that $\Delta^\prime$ is an ISS.  As the link of every vertex of an ISS is a stacked sphere, the same holds in the base space $\Delta.$  But any complex in which the link of every vertex is a stacked sphere is obtained by identifying $\beta_1$ facets of a stacked sphere \cite{Kal}, \cite{Wal}. 
\end{proof}
Since $g_2 \le \binom{g_1+1}{2}$ for any simplicial complex, the above theorem implies that any triangulation of an $S^{n-2}$-bundle over $S^1$ has at least $2n+1$ vertices.  We call a triangulation of an $S^{n-2}$-bundle over $S^1$ with exactly $2n+1$ vertices a {\it minimal identified stacked sphere}, or MISS.  

\section{Constucting ISS's} \label{algorithm}
In this section we show how to construct a MISS for any $n \ge 3.$  Our construction will turn out to be identical to K\"uhnel's minimal triangulations of $S^{2k+1} \times S^1$ and nonorientable $S^{2k}$-bundles over $S^1.$  We also show how with one extra vertex it is possible to triangulate the other sphere-bundles over $S^1.$

Let $\Delta_1$ be the boundary of $\Delta^n,$ the $n$-simplex.  From this point forward we identify the vertices $\{v_1, \dots, v_m\}$ with $\{1, \dots, m\}.$  So, $\Delta_1$ is the complex whose facets are the $n-$subsets of $\{1,\dots,n+1\}.$ Define $\Delta_i$ inductively by setting $\Delta_{i+1}$ equal to the complex obtained by subdividing the facet $\{i+1,\dots, n+i\}$ in $\Delta_i$ with new vertex $n+i+1.$   Evidently each $\Delta_i$ is a stacked sphere.

In order to verify that we can form an ISS we introduce the following notation.  The distance between two vertices $i$ and $j$, denoted $d(i,j)$, is defined to be the minimal length of an edge path between them. To each vertex $i$ we associate the vector 
$$x_i= \left[
  \begin{array}{ c c }
    d(i,1) \\
    . \\
    . \\
    . \\
    d(i,n)
  \end{array} \right],$$
whose entries consist of the distances from vertex $i$ to the vertices $j\in\{1, \dots, n\}$.  In $\Delta_1$ we have initial vectors:
$$x_1 = \left[
  \begin{array}{ c c }
    0 \\
    1 \\
    . \\
    . \\
    . \\
    1
  \end{array} \right], \ 
x_2=\left[
  \begin{array}{ c c }
    1 \\
    0 \\
    . \\
    . \\
    . \\
    1
  \end{array} \right], \dots, 
x_n=\left[
  \begin{array}{ c c }
    1 \\
    1 \\
    . \\
    . \\
    . \\
    0
  \end{array} \right], \ 
x_{n+1} = \left[
  \begin{array}{ c c }
    1 \\
    1 \\
    . \\
    . \\
    . \\
    1
  \end{array} \right] .$$
 
Observe that vertex $n+i$ is introduced in $\Delta_i$  and $x_1, \dots, x_{n+i}$ are unaltered by later subdivisions.  Furthermore, for $i \ge 2, $

\begin{equation}
\label{eqn: min. formula}
 x_{n+i}^{(j)} =  \min \ ( x_{i}^{(j)}, \dots, x_{i+n-1}^{(j)} )+1.
 \end{equation}
 
 \noindent  Here, $x^{(j)}_i$ is the $j$-th coordinate of $x_i.$  Once we reach $\Delta_{2n+1}$ we have the following distance table.
 
 \begin{equation} \label{distance table}
 \begin{array}{ccccccccccc}
 \ &x_1 & \dots & x_{n+1} & x_{n+2} & \dots & x_{2n+1} & x_{2n+2} & x_{2n+3} & \dots & x_{3n+1} \\
 1& 0 & \dots & 1& 2 & \dots & 2 & 3 & 3 & \dots & 3\\
 2& 1 & \dots & 1& 1 & \dots & 2 & 2 & 3 & \dots & 3\\
 \vdots & \vdots & \vdots& \vdots& \vdots& \vdots& \vdots& \vdots& \vdots& \vdots& \vdots \\
 n& 1 & \dots & 1 & 1 & \dots & 2 & 2 & 2 & \dots & 3 \end{array}
 \end{equation}

Examining the last $n$ columns, we see that an ISS can be formed by identifying the pairs of vertices $(1, 2n+2), (2, 2n+3) , \dots, (n, 3n+1).$ 

\begin{defn}
$M^n$ is the MISS obtained by the above identification on $\Delta_{2n+1}.$
\end{defn}

 What are the facets of $M^n$?  At the $i$-th step all $n$-subsets of $\{i,i+1,\dots, n+i\}$ become facets of $\Delta_i,$ except $\{i, i+1, \dots, n+i-1\}$, which has been removed from $\Delta_{i-1}$ to be replaced by the other facets.    However, for every $i, \ 2 \le i \le 2n+1,$ each consecutive subset $\{i,i+1,\dots,i+n-1\}$ is eliminated in $\Delta_{i+1}.$  This leaves $\{1, \dots, n\}$ and $\{2n+2, \dots, 3n+1\}$ as the only consecutive subsets which could be facets.  These are eliminated when the two facets of $\Delta_{2n+1}$ are identified and removed.  Thus, the facets of $M^n$ are all $n$-subsets of the cyclic $2n+1$ sets generated by the set $\{1, \dots, n+1\}$ under the cyclic action  modulo $2n+1,$ other than those which are consecutive.  This is exactly K\"uhnel's generalization of Cs\'asz\'ar's torus.  Hence, when $n$ is odd, $M^n$ is homeomorphic to $S^{n-2} \times S^1,$ and when $n$ is even, $M^n$ is the nonorientable $S^{n-2}$-bundle over $S^1 $ \cite{Ku2}.

Suppose the above algorithm is extended one more step to $\Delta_{2n+2}.$  Then all of the entries of the last two columns are at least three.  So we can form an ISS in many different ways. Two possibilities are to  identify the pairs $(1, 2n+3), \dots, (3n, n-2), (3n+1,n-1), (3n+2, n)$, or use the same pairs, but exchange the last two, i.e. $(1, 2n+3), \dots, (3n, n-2), (3n+1,n), (3n+2, n-1)$.  Since these two identifications have opposite orientations, they must produce the two different $S^{n-2}$-bundles over $S^1. $  As a consequence we have the following.

\begin{thm}
There exist triangulations of $S^{2k} \times S^1$ with $4k+6$ vertices.  There exist triangulations of the nonorientable $S^{2k+1}$-bundle over $S^1$ with $4k+8$ vertices.
\end{thm}

\noindent This was conjectured  by Lutz in \cite{Lu}.  By using identifications which differ by an even permutation it is not hard to see that the $4k+6$ vertex triangulations of $S^{2k} \times S^1$  and the $4k+8$ vertex triangulations of the nonorientable $S^{2k+1}$-bundle over $S^1$ are not unique.

\section{Uniqueness of a MISS} \label{unique}

Throughout this section $\Delta$ is a stacked sphere.  In order to prove that any MISS is combinatorially isomorphic to $M^n$ we require several preliminary results.

\begin{defn}
A {\bf stack} of a stacked sphere is the set of facets created by a maximal subsequence of subdivisions $s_1 \ldots s_k$, where $s_i$ subdivides a facet created by $s_{i-1}$. 
We will sometimes call a stacked sphere with $l$ stacks an $l$-stacked sphere.  The {\bf top} of a stack is the set of facets created by the last subdivision of the stack.  The {\bf top vertex} of the stack is the vertex introduced by $s_k.$
\end{defn}

\begin{figure} 
 \scalebox{0.55}[0.55]{\includegraphics{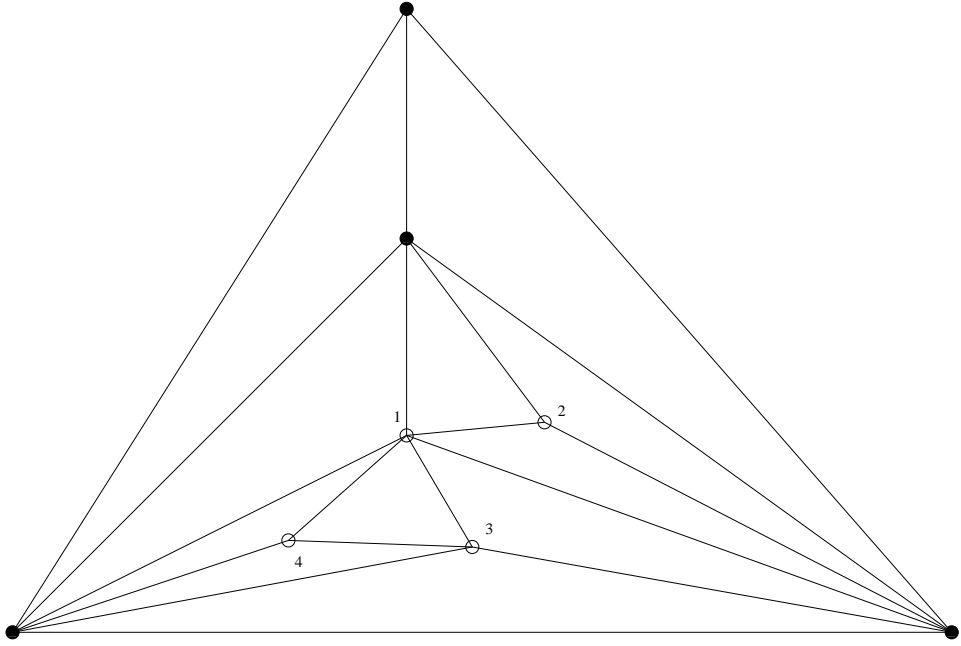}}
  \caption{A $2$-stacked sphere} \label{fig1}
\end{figure}

\begin{example}
The two-dimensional sphere in Figure \ref{fig1} is a $2$-stacked sphere.  The unlabeled vertices are the original vertices of the  boundary of the tetrahedron. One stack has top vertex $2$  and is formed by the subdivisions corresponding to vertices $1$ and $2.$  The other stack is formed by the subdivisions for vertices $1$, $3$, and $4$ and its top  vertex is $4.$
\end{example}

\begin{prop} \label{at most two}
Suppose a MISS is obtained from a stacked sphere $\Delta$ by handle addition.  Then $\Delta$ has at most two stacks. Further, the top of every stack must have one of the two identified facets.
\end{prop}

\begin{proof}
At least one facet from the top of each stack must be identified.  Otherwise, the vertex corresponding to the last subdivision in the stack could be removed from the MISS, contrary to its definition.  Since the tops of distinct stacks are disjoint, there can be at most two stacks.
\end{proof}

\begin{prop} \label{1-stack}
If $M$ is a \miss, then $M$ can be obtained by handle addition on $\Delta,$ where $\Delta$ has only one stack.  
\end{prop}

\begin{proof}
By Theorem \ref{MISS->SS} and Proposition \ref{at most two}  there exists $\Delta^\prime$ a stacked sphere with one or two stacks so that $M$ is obtained by handle addition on $\Delta^\prime$. If $\Delta^\prime$ has one stack we are done, so assume that $\Delta^\prime$ has two stacks. 

Denote the vertices of the identified facets by $\{1,\dots,n\}$ and $\{1^\prime, \dots, n^\prime\},$ where vertex $i$ is identified with vertex $i^\prime$ and $1$ is the top vertex of its stack.  Designate the remaining vertices at the top of their respective stacks, $n+1$ and $(n+1)^\prime.$  Now, we undo the subdivision that created $1$ and subdivide the facet  $\{1', \ldots, n'\}$ with new vertex $(n+1)^{\prime\prime}$.  Call this new stacked sphere $\Delta^{\prime\prime}.$  For $2 \le i \le n$ the distance between $i$ and $i^\prime$ is still at least three.  In addition, $d(n+1, (n+1)^{\prime\prime})$ is at least three.  If not, then there is an edge between $n+1$ and some $i^\prime,$ which is impossible as that would imply that $d(i,i^\prime) \le 2.$  It is now easy to check that the MISS obtained by identifying $i$ with $i^\prime$ for $2 \le i \le n$, and  the pair of vertices $(n+1), (n+1)^{\prime\prime}$ is combinatorially isomorphic to $M$.  Evidently we can repeat this procedure until one of the stacks is gone, resulting in the desired $\Delta.$
 \end{proof}
 
 \begin{thm} \label{uniqueness}
 Let $M$ be an $(n-1)$-dimensional MISS.  Then $M$ is isomorphic to $M^n.$
 \end{thm}
 
 \begin{proof}
 By Theorem \ref{MISS->SS} and Proposition \ref{1-stack} we can assume that $M$ is obtained by handle addition from a 1-stacked sphere $\Delta.$ Label the vertices of $\Delta$ with $\{1,\dots,3n+1\},$ where $\{1,\dots,n+1\}$ is the original boundary of the $(n-1)$ simplex and the order of the vertices reflects the order of the subdivisions.  Without loss of generality we can assume that the first subdivided facet is $\{2,\dots,n+1\}.$  Since $\Delta$ is a one-stacked sphere no other vertex will share an edge with $1.$  Hence, one of the identified facets must contain the vertex 1.  Otherwise, the closed star of 1 could be removed from $M$ and replaced by $\{2, \dots, n+1\},$ contradicting the fact that $M$ is a MISS. Thus one of the identified facets consists of $n$ of the first $n+1$ vertices (including the vertex 1), and the other contains the top vertex of the stack, vertex $3n+1.$  
 
    Now consider the distances $d(i,j)$ for $1 \le i \le n+1$ as in (\ref{distance table}), but with an extra row for vertex $n+1.$  Since $\Delta$ is a 1-stacked sphere the table can be constructed as follows.  It begins
    
    $$\begin{array}{ccccccc}
     & x_1 & x_2 & \dots & x_n & x_{n+1} & x_{n+2} \\
    1& 0 & 1 & \dots & 1 & 1 & 2\\
    2& 1 & 0 & \dots & 1 & 1 & 1\\
    3& 1 & 1 &\dots & 1 & 1& 1 \\
\vdots&\vdots&\vdots&\vdots&\vdots&\vdots&\vdots \\
    n & 1 & 1 & \dots & 0 & 1 & 1\\
   n+1 & 1 & 1 & \dots & 1 & 0 & 1.
    \end{array}
    $$
 \noindent  Since vertex $n+2$ was formed by subdividing the facet $\{2,\dots,n+1\},$ the $x_{n+2}$ column is obtained by ``crossing off" the first column and applying (\ref{eqn: min. formula}) to the remaining $n$ columns.  Since $\Delta$ is a 1-stacked sphere, each successive column is obtained in the same fashion.  Delete one previous column and apply the analog of (\ref{eqn: min. formula}) to the remaining $n$ columns.  The crossed off column corresponds to the vertex in the current top of the stack which is not in the newly subdivided facet.  
 
We claim that the first $n$ columns crossed off are a subset of those headed by $x_1, \dots, x_{n+1}.$ Suppose not.  In view of (\ref{eqn: min. formula}) this implies there are two indices $i$ and $j$ such that $2 \le i<j \le n+1$ and that the $i^{th}$ and $j^{th}$ coordinates of $x_{2n+1}$ and all previous $x_k$ are $1$ or $0$.  However, this makes it impossible for any of the remaining vectors to have a 3 in either of the $i^{th}$ or $j^{th}$ coordinates.  In particular, no vertex may be identified with vertices $i$ or $j,$ contrary to the fact that one of the identified facets contains vertex 1.  So, renumbering the vertices $\{2,\dots,n+1\}$  if necessary, the first $n$ subdivisions are identical to those made in constructing $\Delta_{n+1}$ and the distance table for $\Delta$ begins
 
 \begin{equation}
 \begin{array}{ccccccc}
 \ &x_1 & \dots & x_{n+1} & x_{n+2} & \dots & x_{2n+1} \\
 1& 0 & \dots & 1& 2 & \dots & 2 \\
 2& 1 & \dots & 1& 1 & \dots & 2  \\
 \vdots & \vdots & \vdots& \vdots& \vdots& \vdots& \vdots \\
 n& 1 & \dots & 1 & 1 & \dots & 2 \\
 n+1 & 1 & \dots & 0 & 1 & \dots & 1
  \end{array}
 \end{equation}  
 
 At the next step the column headed by $x_{n+1}$ must be crossed off.  If not, then $x_{2n+2}$ and all of its predecessors do not contain any coordinate with a value of $3$ or greater.  This would leave only $n-1$ vertices which could be identified.  Furthermore, since $x_{2n+2} = (3,2,\dots,2),$ vertex $2n+2$ must be identified with vertex 1.  Now repeating the argument in each of the subsequent columns shows that $\Delta$ and $M$ are constructed in exactly the same fashion as $\Delta_{2n+1}$ and $M^n.$
 \end{proof}
 
 By Theorem \ref{MISS->SS} and the argument in Corollary \ref{cor}  we immediately obtain the following result.
 
 \begin{thm}
If $n \ge 5,$ then  $M^n$ is the unique triangulation of any $(n-1)$-dimensional homology manifold with $2n+1$ vertices, $ \beta_1 \neq 0$ and whose orientation double cover has vanishing second Betti number.  
 \end{thm}
 
 \section{All possible $(f_0,f_1)$ pairs} \label{edges}
 
 An $S^1$-bundle over $S^1$ is a torus or Klein bottle, and their possible $f$-vectors are well known.  Walkup provided a characterization of all $f$-vectors of $S^2$-bundles over $S^1$ \cite{Wal}.
 Theorems \ref{MISS->SS} and \ref{uniqueness} determine the minimum number of vertices and edges for any $S^k$-bundle over $S^1$ when $k \ge 3.$ All possible pairs $(f_0,f_1)$ for triangulations of $S^{2k+1} \times S^1$  were established in \cite{Sw10}.  To do this for the other spherical bundles over $S^1$  we use sequences of bistellar moves. 

 Suppose $\Delta$ is a triangulation of an $(n-1)$-manifold and let $A$ and $B$ be disjoint subsets of vertices of  $\Delta$ with $|A|=2$ and $ |B|=n-1.$ Furthermore, assume that the induced subcomplex on $A \cup B$ is the suspension of $B.$  Equivalently, a subset of $A \cup B$ is a face of $\Delta$ if and only if it does not contain both vertices in $A.$  In this case the induced subcomplex is an $(n-1)$-ball with boundary $\partial A \ast \partial B,$ so we can replace it with $A \ast \partial B,$ an $(n-1)$-ball with the same boundary, without changing the homeomorphism type of $\Delta.$  As long as $n \ge 4$ the vertices of the new complex will be the same, and the edge set will be the edge set of $\Delta$ with one new edge between the vertices of $A.$ This transformation is called a {\it bistellar move} on $A \cup B.$  
  
\begin{thm}  \label{edges+}
Let $\Delta$ be a triangulation of an $S^k$-bundle over $S^1, k \ge 2.$  If $k$ is odd and $\Delta$ is orientable, or $k$ is even and $\Delta$ is nonorientable, then $f_0 \ge 2k+5.$  Otherwise, $f_0 \ge 2k+6.$  Both bounds are sharp and in all  cases, if there exists a triangulation with $f_0$ vertices, then there exists a triangulation with $f_1$ edges and $f_0$ vertices if and only if $(k+2) f_0 \le f_1 \le \binom{f_0}{2}.$

\end{thm}

\begin{proof}
  For $k=2,$ see \cite{Wal}. So assume for the rest of the proof that $k \ge 3.$ The  minimum values for $f_0$ follow from Corollary \ref{cor} and Theorem \ref{uniqueness}.  The  lower bounds for $f_1$ are given by Corollary \ref{cor}, while the upper bond holds trivially for any simplicial complex.  It remains to show that there exists a triangulation for each such $f_1.$

Set $n=k+2,$ so $n \ge 5$ and is the cardinality of every facet.  Suppose $f_0$ satisfies the above inequalities. First we can construct an ISS homeomorphic to the desired bundle by taking the stacked sphere $\Delta_{f_0}$ and identifying the pairs of vertices $\{1,f_0+1\}, \{2, f_0 + 2\}, \dots, \{n-1, f_0+n-1\},\{n, f_0+n\}$ or the pairs of vertices $\{1,f_0+1\}, \{2, f_0 + 2\}, \dots, \{n-1, f_0+n\},\{n, f_0+n-1\},$ where the choice of pairing depends on whether the orientable or nonorientable bundle is under consideration.  Since any stacked sphere with $f_0+n$ vertices has $\binom{n+1}{2} + (f_0-1)n$ edges, the resulting ISS will have $n f_0$ edges.  Thus, it is suffices to show how to perform successive bistellar moves on this ISS until the $1$-skeleton is the complete graph on the $f_0$ vertices.  For the first pairing this was done in \cite[Theorem 5.2]{Sw10}, but we repeat a similar argument here for the convenience of the reader.

For the first pairing the missing edges in the resulting MISS consist of all pairs of vertices $\{i,j\}$ with $n< j-i < f_0 - n.$  Equivalently, the cyclic distance between $i$ and $j$ is at least $n+1.$  Group the pairs of nonedges according to the size of $j-i.$  So the first group consists of $\{1,n+2\}, \{2,n+3\}, \dots,$ the second group $\{1,n+3\}, \{2,n+4\},$ etc.  The first bistellar move uses $A = \{1, n+2\}$ and $B = \{2,3,5,6,\dots, n+1\}.$ The next bistellar move uses $A=\{2,n+3\}$ and $B = \{3,4,6,7,\dots, n+2\}.$  We continue with bistellar moves which introduce an edge between each pair  of nonedges in the first group with $A$ the pair $\{i,j\}$ and $B$ consisting of $\{i+1,i+2\}$ and the $n-3$ vertices preceding  $j.$   

In order to get edges for the second group of nonedge pairs we proceed as follows.  Let $A=\{1,n+3\}$ and $B=\{2,3,6,7,\dots,n+2\}.$ The previous bistellar move using $\{2,n+3\}$ and $\{3,4,6,7,\dots,n+2\}$ makes it possible to  execute a bistellar move with this $A$ and $B.$ Similarly, it is possible to use $A=\{2,n+4\}$ and $B=\{3,4,7,8,\dots,n+3\}$ to force an edge between $2$ and $n+4$ because of the previous bistellar move on $\{3,n+4\}$ and $\{4,5,7,8,\dots,n+3\}.$ Each pair $\{i,j\}$ in the second group uses the bistellar move associated to $\{i+1,j\}$ in the first group.  This procedure allows us to take care of all of the nonedge pairs in the second group.  Now each succeeding group of nonedges repeats the process of using the proceeding group's bistellar moves in order to perform their own bistellar moves until there are no more nonedge pairs remaining and the 1-skeleton is the complete graph.

What are the nonedges in the ISS for the second pairing?  With two exceptions they are the same as before.  The vertex $n$, which was previously identified with $f_0+n,$ is now paired with $f_0+n-1,$ and hence has an edge with $f_0-1$ which is not true for the first pairing.  Similarly, vertex $n-1$ is now identified with $f_0 +n$ instead of $f_0+n-1,$ so it does not have an edge with $f_0-1,$ which it did in the first identification scheme. The nonedge pair $\{n-1,f_0-1\}$ lies in the last group.  Since there is no bistellar move corresponding to $\{n,f_0-1\}$ (which would have been in the next to last group), a different bistellar move is required to get an edge between vertex $n$ and $f_0.$  Otherwise, we use exactly the same algorithm as above.  We connect the vertices $n-1$ and $f_0-1$ in the last bistellar move using $A = \{n-1,f_0-1\}$ and $B = \{f_0,f_0+1, f_0+3, \dots,f_0+n-2, f_0+n-1\} = \{f_0,1,3,\dots,n-2,n\}.$

\end{proof}

\begin{remark}
The same argument as in \cite[Theorem 4.7]{Sw10} shows that for a fixed $S^k$-bundle over $S^1$ any ISS which minimizes $f_0$ and $f_1$ has the minimum number of faces in every dimension over all triangulations of the bundle.  
\end{remark}

Since the $f$-vector of a 4-manifold is completely determined by $g_1$ and $g_2,$ the above theorem immediately  gives a complete list of all possible $f$-vectors of the nonorientable $S^3$-bundle over $S^1.$  

\begin{thm} \label{s3 twist s1}
The following are equivalent.
\begin{enumerate}
  \item
  $(g_0,g_1,g_2)$ is the $g$-vector of a triangulation of the nonorientable $S^3$-bundle over $S^1.$
  \item
  $g_0 = 1, g_1 \ge 6, 15 \le g_2 \le \binom{g_1+1}{2}.$
  \end{enumerate}
\end{thm}
    
{\bf Acknowledgment.} All three authors were partially supported by Cornell's 2006 summer REU program.  The last author was also partially supported by DMS-0600502.


\begin{thebibliography}{10}

\bibitem{BD}
B.~Bagchi and Basudeb D.
\newblock The lower bound theorem and minimal triangulations of sphere bundles
  over the circle, 2006.
\newblock arXiv:math.GT/0610829.

\bibitem{BKu}
U.~Brehm and W.~{K\"uhnel}.
\newblock Combinatorial manifolds with few vertices.
\newblock {\em Topology}, 26(4):465--473, 1987.

\bibitem{JR}
M.~Jungerman and G.~Ringel.
\newblock Minimal triangulations of orientable surfaces.
\newblock {\em Acta Math.}, 145:121--154, 1980.

\bibitem{Kal}
G.~Kalai.
\newblock The diameter of graphs of convex polytopes and f-vector theory.
\newblock In {\em Applied geometry and discrete mathematics}, pages 387--441.
  Amer. Math. Soc., Providence, RI, 1991.
\newblock DIMACS Ser. Disc. Math. Theoret. Comput. Sci., 4.

\bibitem{Ke}
V.~Klee.
\newblock A combinatorial analogue of {Poincare's} duality theorem.
\newblock {\em Canadian J. Math.}, 16:517--531, 1964.

\bibitem{Ku2}
W.~K{\"u}hnel.
\newblock Higher dimensional analogues of {Cs\'asz\'ar's} torus.
\newblock {\em Result. Math.}, 9:95--106, 1986.

\bibitem{Lu}
F.~Lutz.
\newblock Triangulated manifolds with few vertices: {C}ombinatorial manifolds,
  2005.
\newblock arXive: math.CO/0506372.

\bibitem{Ri}
G.~Ringel.
\newblock Wie man die geschlossenen nichtorientierbar {Fl\"achen} in
  {m\"oglichst} wenig {D}reiecke zerlegen kann.
\newblock {\em Math. Ann.}, 130:317--326, 1955.

\bibitem{Ste}
N.~E. Steenrod.
\newblock The classification of sphere bundles.
\newblock {\em Ann. of Math.}, 45:294--311, 1944.

\bibitem{Sw10}
E.~Swartz.
\newblock Face enumeration: from spheres to manifolds, 2005.
\newblock \url{http://www.math.cornell.edu/~ebs/spherestomanifolds.pdf}.

\bibitem{Wal}
D.~Walkup.
\newblock The lower bound conjecture for 3 and 4 manifolds.
\newblock {\em Acta Math.}, 125:75--107, 1970.

\end{thebibliography}
\end{document}